\documentstyle[11pt]{article}
\newtheorem {theorem}{Theorem}[section]

\newtheorem {lemma}[theorem]{Lemma}
\newtheorem {corollary}[theorem]{Corollary}
\newcounter{conjecture}
\newcounter{remark}

\newcommand{\eqnsection}{
    \renewcommand{\theequation}{\thesection.\arabic{equation}}
    \makeatletter
    \csname @addtoreset\endcsname{equation}{section}
    \makeatother}

\def \be{\begin{equation}}
\def \ee{\end{equation}}
\def \bt{\begin{theorem}}
\def \et{\end{theorem}}
\def \bea{\begin{eqnarray}}
\def \eea{\end{eqnarray}}
\def \bas{\begin{eqnarray*}}
\def \eas{\end{eqnarray*}}

\def \ga{\gamma}
\def \Ga{\Gamma}
\def \de{\delta}

\def \ep{\epsilon}

\def \la{\lambda}
\def \La{\Lambda}

\def \th{\theta}

\def \ze{\zeta}

\def \ff{\infty}

\def \wt{\widetilde}

\def \rar{\rightarrow}

\def \CC{{\cal C}}

\def \HH{{\cal H}}

\def \LL{{\cal L}}

\def \({\left(}
\def \){\right)}

\def \lc{\left\{}
\def \rc{\right\}}
\def \bsq{\hfil $\Box$}
\def \nn{\nonumber}

\def \bc{\begin{center} }
\def \ec{\end{center} }
\begin{document}

\eqnsection
\newcommand{\Ini}{{I_{n,i}}}
\newcommand{\reals}{{I\!\!R}}
\newcommand{\complex}{{I\!\!\!C}}
\newcommand{\F}{{\cal F}}
\newcommand{\D}{{\cal D}}
\newcommand{\Fn}{{{\cal F}_n}}
\newcommand{\Hn}{{{\cal H}_n}}
\newcommand{\Fp}{{{\cal F}^p}}
\newcommand{\PPP}{{\bf P}}
\newcommand{\Pop}{{P\otimes \PPP}}
\newcommand{\hm}{\HH^\varphi}
\newcommand{\nuw}{{\nu^W}}
\newcommand{\ths}{{\theta^*}}
\newcommand{\beq}{\begin{equation}}
\newcommand{\eeq}{\end{equation}}
\newcommand{\integers}{{\rm I\!N}}
\newcommand{\E}{{\rm I\!E}}
\newcommand{\te}{{\tilde{\delta}}}
\newcommand{\tI}{{\tilde{I}}}
\newcommand{\loge}{{\log(1/\ep)}}
\newcommand{\logen}{{\log(1/\ep_n)}}
\newcommand{\epn}{{\ep_n}}
\def\var{{\rm Var}}
\def\cov{{\rm Cov}}
\def\one{{\bf 1}}
\def\leb{{\cal L}eb}
\def\Ho{{\mbox{\sf H\"older}}}
\def\thi{{\mbox{\sf Thick}}}
\def\thin{{\mbox{\sf Thin}}}
\def\dimm{{\overline{{\rm dim}}_{_{\rm M}}}}
\def\dimp{\dim_{_{\rm P}}}
\def\psih{\psi}
\def \MU{{\cal J}}
\bibliographystyle{amsplain}

\title{{\bf \LARGE  Frequently visited sets for random walks }}

\author{Endre Cs\'{a}ki\thanks{Research supported by the Hungarian
National Foundation for Scientific Research, Grant No. T 037886 and  T
043037.}\,\,\, Ant\'{o}nia  F\"{o}ldes\thanks{Research supported by a
PSC CUNY Grant,  No. 65685-0034.}\,
  \, \,P\'{a}l R\'{e}v\'{e}sz$^*$\, \, Jay Rosen\thanks {Research
supported, in part, by grants from the NSF  and from PSC-CUNY.} \, \,
Zhan Shi }
\date{}

\maketitle

\begin{abstract} We study the occupation measure of various sets for  a
symmetric transient random walk in $Z^d$ with finite variances.
Let $\mu^X_n(A)$ denote the occupation time of the set $A$ up to
time $n$. It is shown that $\sup_{x\in Z^d}\mu_n^X(x+A)/\log n$
tends to a finite limit as $n\to\infty$. The limit is expressed in
terms of the largest eigenvalue of a matrix involving the Green's
function of $X$ restricted to the set $A$. Some examples are
discussed and the connection to similar results for Brownian
motion is given.
\end{abstract}



\section{Introduction}\label{secintro}

Let $X_{ n},\,n=0,1,\ldots$ be a symmetric transient random walk in
$Z^d$ ($d\geq 3$). We will always assume that $X_{ n},\,n=0,1,\ldots$ is
not supported on any subgroup strictly smaller than $Z^d$. We denote by
$\mu^X_n$  its {\em occupation measure}:
\[
\mu^X_n(A)=\sum_{ j=0}^n \one_A(X_{ j})
\] for all  sets $A\subseteq Z^d$.  Let
$q_{ n}(x)=\PPP(X_{ n}=x )$.  As usual, we let
\begin{equation} G(x)=\sum_{ k=0}^{ \ff}q_{ k}(x)\label{10.1}
\end{equation} denote the Green's function for $\{ X_n \}$.  For any finite
$A\subseteq Z^d$ let $\La_{A}$ denote the largest eigenvalue of the
$|A|\times |A|$ matrix
\be G_{ A}(x,y)=G(x-y),\hspace{ .2in}x,y\in A.\label{matrix}\ee

\bt\label{theo-basic} If $X$ has finite second moments then \be
\lim_{n \to \ff} \sup_{x\in Z^{ d}} {\mu^X_n(x+A) \over \log n}
=-1/\log (1-1/\La_{A}) \hspace{.2in} {\rm a.s.} \label{seq-amir3}
\ee and \be \lim_{n \to \ff}  \sup_{0\leq m \leq n} {\mu^X_n(X_{
m}+A) \over \log n} = -1/\log (1-1/\La_{A}) \hspace{.2in} {\rm
a.s.} \label{eq-slevy} \ee \et

  For our first example, when $A=\{ 0\}$, $\La_{A}=G( 0)=1/\ga_{ d}$,
where
$\ga_{ d}$
  is the probability of no-return to the origin, and in the case of  the
simple random walk we recover Theorem 13 of \cite{ET}.

Here are some other examples. Set $t_{ y}=\PPP(T_{ y}<\ff)$, where
$T_{y} := \inf \{s > 0: X_s=y\}$.  Let $S( 0,1)=\{ e_{1},\ldots, e_{ d},-e_{
1},\ldots,-e_{ d}\}$, $B( 0,1)=\{0\}\cup S( 0,1)$, be the (Euclidean)
sphere and ball in
$Z^{ d}$ of radius $1$ centered at the origin.

\bt\label{theo-examp} If $X$ has finite second moments, then for any
$0\neq y\in Z^{ d}$
\be
\lim_{n \to \ff} \sup_{x\in Z^{ d}} {\mu^X_n(x+\{ 0,y\}) \over
\log n} =-1/\log (1-\ga_{ d}/( 1+t_{ y}))
\hspace{.2in} {\rm a.s.}
\label{ex.1}
\ee

For the simple random walk
\be
\lim_{n \to \ff} \sup_{x\in Z^{ d}} {\mu^X_n(x+S( 0,1)) \over
\log n} =-1/\log \(1-\ga_{ d}/2d( 1-\ga_{ d})\)
\hspace{.2in} {\rm a.s.}
\label{ex.2}
\ee and
\be
\lim_{n \to \ff} \sup_{x\in Z^{ d}} {\mu^X_n(x+B( 0,1)) \over
\log n} =-1/\log \({p+\sqrt{p^{ 2}+2/d}
\over 2}\)
\hspace{.2in} {\rm a.s.}
\label{ex.3}
\ee where $p=1-1/2d( 1-\ga_{d})$.
\et

\begin{corollary}
\label{less}
If $X$ has finite second moments, then for any fixed $K>0$
\[
\lim_{n\to\infty}\frac{\max_{x,y\in Z^d: |x-y|\leq
K}(\mu_n^X(\{x,y\})}{\log n}<-2/\log(1-\gamma_d)
\hspace{1cm} \mathrm{a.s.}
\]
\end{corollary}

Since the constant for one-point set in Theorem 1.1 is
$-1/\log(1-\gamma_d)$, this corollary expresses the fact that any two
points with individual occupation measures up to time $n$, both close to
the maximum, should be at a distance larger than any constant $K>0$.
In particular, a neighbor of a maximally visited point is not maximally
visited.

Let $W_t$ denote Brownian motion in $\reals^d$, $d\geq 3$. We denote
by
$\nu^{ W}_T$
   its {\em occupation measure}:
\[
\nu^{ W}_T(A)=\int _0^T \one_A(W_t)\,dt
\] for all Borel sets $A\subseteq \reals^d$.  Let $K\subseteq
\reals^d$ be a fixed compact neighborhood of the origin which is the
closure of its interior and set
$K(x,r)=x+rK$.

As usual, we let
\begin{equation} u^0(x)={c_{d}\over |x|^{d-2}}\label{10.2}
\end{equation} denote the $0$-potential density for $\{ W_t \}$, where
$c_{d}=2^{-1} \pi^{-d/2}
\Ga({d\over 2}-1)$. Let $\La^{ 0}_K$ denote the norm of
\[ R_{K}f(x)=\int_{K} u^0(x-y)f(y)\,dy
\] considered as an operator from $L^2\(K,\,dx\)$ to itself. If
$B(x,r)$ denotes the Euclidean ball in $R^d$ of radius $r$ centered  at
$x$, it is known, \cite{Ciesielski-Taylor}, that
$\La^{ 0}_{B(0,1) }=2r^{- 2}_{ d}$ where
$r_{ d}$ is the smallest positive root of the Bessel function $J_{ d/2-2}$.

We mention that it can be shown, at least for $K$ convex, that for any $S
\in (0,\infty)$ and any $T \in (0,\infty]$,
\be
\lim_{\ep \to 0} \sup_{|x| \leq S} {\nu_T(K(x,\ep)) \over
\ep^2|\log \ep|} = 2 \La^{ 0}_{K} \hspace{.2in} {\rm a.s.}
\label{seq-amir3j}
\ee and
\be
\lim_{\ep \to 0} \sup_{0\leq t \leq T} {\nu_T(K(W_{ t},\ep)) \over
\ep^2|\log
\ep|}
  = 2 \La^{ 0}_{K}
\hspace{.2in} {\rm a.s.}
\label{eq-slevyj}
\ee (These results are mentioned for motivation. They are not used in  the
rest of the paper).

For any $x\in R^{ d}$ and $\ep>0$, let
$e_{ \ep}( x)=x+[0,\ep]^{ d}$, the cube of edgelength $\ep$ with  `lower'
corner at
$x$. Set
\begin{equation}
\LL_{ \ep}( K) =\{ x \in \ep Z^{ d}\,|\,  e_{ \ep}( x)\subseteq K  \},
\hspace{ .2in}\mbox{ and }\hspace{ .2in}
\CC_{ \ep}( K)=\bigcup_{x\in  \LL_{ \ep}( K)}e_{ \ep}( x)\label{1.edgea}
\end{equation} and assume that
\begin{equation}
\lim_{ \ep\rar 0}\la^{ d}(\CC_{ \ep}( K) )=\la^{ d}(K )\label{1.edgeb}
\end{equation} where $\la^{ d}$ denotes Lebesgue measure.

Note that $\ep^{ -1}\LL_{\ep}( K)\subseteq Z^{ d}$.

\bt\label{theo-conv} Assume that  $X_{ 1}$ has $d-1$ moments
  and covariance matrix equal to the identity. Then
\label{theo-approx}
\be
\lim_{\ep \to 0} \ep^{ 2}\La_{\ep^{ -1}\LL_{\ep}( K)}  =  \La^{ 0}_{K}
\label{1.edgec}
\ee and consequently
\begin{equation} -\lim_{\ep \to 0} \ep^{ 2}/\log (1-1/\La_{\ep^{
-1}\LL_{\ep}( K)}) =  \La^{ 0}_{K}.\label{1.edgez}
\end{equation}
\et

Section \ref{seclocal} states and proves the  crucial Localization
Lemma \ref{lb2}. Theorem \ref{theo-basic} is proven in section
\ref{sec-1.1}, Theorem \ref{theo-examp} and Corollary \ref{less}
are proven in section \ref{sec-ex}, and Theorem \ref{theo-conv} is
proven in section \ref{sec-1.2}.

\section{Localization for random walk occupation
measures}\label{seclocal}

We start by providing a convenient representation of the law of the total
occupation measure $\mu^X_{\infty}\(A\)$. This representation is the
counterpart of the Ciesielski-Taylor representation for the total occupation
measure of spatial Brownian motion in
\cite[Theorem 1]{Ciesielski-Taylor}.

Let $( f,g)_{ A}=\sum_{ x\in A}f( x)g( x)$, and let $\de_{ 0}( x)$ be  the
function on $A$ defined by $\de_{ 0}( x)=\de(0, x),\,x\in A$.

\begin{lemma}\label{repres} Let $\{ X_n \}$ be a symmetric transient
random walk in
$Z^d$, and let $A$ be a finite set in $Z^{ d}$ which contains the origin.
Then,

\be
\PPP\(\mu^X_\infty (A) > u \)=\sum_{j}h_j
\( {\la_{ j}-1 \over \la_{ j}}\)^{u}\hspace{ .2in}u=0,1,\ldots,
\label{new:rep}
\ee  where $\la_1>\la_2 \geq \cdots \geq \la_{ |A|}\geq 1/2$ are the
eigenvalues of the symmetric matrix $G_{ A}$ with the corresponding
orthonormal eigenvectors
$\phi_j(y)$, $h_j:=  (1,\phi_j)_{A}\phi_j(0)$.
\end{lemma}

\noindent {\bf Proof of Lemma \ref{repres}:} Let $\MU =\mu^X_\infty
(A)$ and set
$\bar G(x-y)=\sum_{ k=1}^{ \ff}q_{ k}( x-y)=G(x-y)-q_{ 0}( x-y)$.

Note that for any $m $,
\bea
\E \( \MU^m\) &=&
\E \(\lc \sum_{ i=0}^{\ff} \one_{A}(X_i)\rc^m\) =\sum_{ i_{ 1},\ldots,
i_{m}=0}^{ \ff}\E \( \prod_{ j=1}^{m} \one_{A}(X_{i_{ j} })\)\nn\\ &=&
\sum_{ k=1}^{ m}\sum_{\stackrel{c_{ 1},\ldots, c_{ k}\in [1,m]}{c_{
1}+\cdots+ c_{ k}=m}}{m \choose c_{ 1},\ldots, c_{
k}}\sum_{A^k}\sum_{0\leq n_1<\cdots< n_k<\ff}\,\,
\prod_{j=1}^k q_{n_j-n_{j-1}}(x_j-x_{j-1}).
\label{11.5}
\eea Here, $k$ is the number of distinct indices $n_1<\cdots<n_k$ among
the indices
$i_{ 1},\ldots, i_{m}$ and $c_{l}$ is the number of times that
$n_{ l}$ appears, i.e. $c_{l}=\#\{ j;\,1\leq j\leq m,
\,\,i_{j}=n_{ l}\}
$. The factor ${m \choose c_{ 1},\ldots, c_{ k}}$ is the number of
ways to assign  the value $n_{ l}$ to $c_{l}$ of the indices $i_{
1},\ldots, i_{m}$, for each $1\leq l\leq k$.

Also, we have that
\begin{eqnarray}
\sum_{A^k}\sum_{0\leq n_1<\cdots<n_k<\ff}\,\,
\prod_{j=1}^k q_{n_j-n_{j-1}}(x_j-x_{j-1})  &=&\sum_{A^k}G(x_1)
\prod_{j=2}^k \bar G(x_j-x_{j-1})\label{11.5add}\\ &=& (1,\bar G_{ A}^{
k-1}G_{ A}\de_{ 0})_{A}.\nonumber
\end{eqnarray}
   Hence (we justify the computations shortly)
\begin{eqnarray}
\E \( e^{\ze \MU }\) &=& 1+\sum_{ m=1}^{ \ff}{ \ze^{ m}\over m!}
\sum_{ k=1}^{ m}\sum_{\stackrel{c_{ 1},\ldots, c_{ k}\in [1,m]}{c_{
1}+\cdots+ c_{ k}=m}}{m \choose c_{ 1},\ldots, c_{ k}}(1,\bar G_{ A}^{
k-1}G_{ A}\de_{ 0})_{A}\label{11.5s}\\  &=& 1+\sum_{ k=1}^{ \ff}
\sum_{ m=k}^{ \ff}\sum_{\stackrel{c_{ 1},\ldots, c_{ k}\in [1,m]}{c_{
1}+\cdots+ c_{ k}=m}}\prod_{ j=1}^{ k}{\ze^{ c_{ j}} \over c_{ j}!}(1,\bar
G_{  A}^{ k-1}G_{ A}\de_{ 0})_{A}
\nonumber\\  &=& 1+\sum_{k=1}^{ \ff} (e^{\ze }-1)^{ k}(1,\bar G_{ A}^{
k-1}G_{ A}\de_{ 0})_{A}
\nonumber
\end{eqnarray}

  $G_{ A}$ is a symmetric matrix. Let $\psi( p)$ denote the  characteristic
function of $X_{ 1}$. Then $\psi( p)$ is real and $|\psi( p)|\leq 1$. Thus
$0\leq 1-\psi( p)\leq 2$, or equivalently ${ 1\over 1-\psi( p)}\geq {1
\over  2}$. Hence, using the Fourier transform representation
$G(x-y)= \int e^{i((x-y)\cdot p)}( 1-\psi( p))^{-1}\,dp$ we can see that
$\sum_{ x,y\in A}G_{ A}( x,y)a_{ x}a_{ y}\geq {1 \over 2}\sum_{ x\in
A}a^{ 2}_{ x}$ for any $\{a_{ x}\in R^{ 1}; x\in A \}$.
   By the standard theory for symmetric matrices,
$G_{ A}$ has  all eigenvalues $\geq 1/2$, and the corresponding
eigenvectors of
$G_{ A}$, denoted $\{ \phi_j
\}$ form an orthonormal basis of
$L^2\(A\)$ (see \cite[Theorems VI.15, VI.16]{Reed-Simon2}). Moreover,
since the entries of $G_{ A}$ are strictly positive,
  by the Perron-Frobenius Theorem, see
\cite[Theorem XIII.43]{Reed-Simon}, the eigenspace corresponding to
$\La_{ A}=\la_1$ is one dimensional, and we may and shall choose
$\phi_1$
  such that $\phi_1(y)>0$ for all $y \in A$.

Thus we can write (\ref{11.5s}) as
\begin{eqnarray}
\E \( e^{\ze \MU }\) &=& 1+\sum_{k=1}^{ \ff} (e^{\ze }-1)^{ k}(1,\bar
G_{ A}^{ k-1}G_{ A}\de_{ 0})_{A}\label{11.5t}\\  &=& 1+ \sum_{j=1}^{
|A|}(1,\phi_j)_{A}(\phi_j,\de_{ 0})_{A}\sum_{k=1}^{ \ff} (e^{\ze  }-1)^{
k} ( \la_{ j}-1)^{ k-1}\la_{ j}
\nonumber\\  &=& 1+ \sum_{j=1}^{ |A|}h_j (e^{\ze }-1)\la_{
j}\sum_{k=1}^{ \ff} (e^{\ze }-1)^{ k-1} ( \la_{ j}-1)^{ k-1}
\nonumber
\end{eqnarray} where $h_j=(1,\phi_j)_{A}(\phi_j,\de_{ 0})_{A}$.  It is
now easy to see that we can justify the derivation of (\ref{11.5s}) and
(\ref{11.5t}) if
\begin{equation} |(e^{\ze }-1) ( \la_{j}-1)|<1,\hspace{ .2in}\forall
j.\label{test1}
\end{equation} In that case we can write (\ref{11.5t}) as
\begin{equation}
\E \( e^{\ze \MU }\)=1+
\sum_{j}h_j{(e^{\ze }-1)\la_{ j} \over 1-(e^{\ze }-1) ( \la_{
j}-1)}.\label{11.5uj}
\end{equation}

Since $\sum_{j=1}^{ |A|}h_{ j}=\sum_{j=1}^{
|A|}(1,\phi_j)_{A}(\phi_j,\de_{ 0})_{A}=(1,\de_{0})_{A}=1$  we have that

\begin{equation}
\E \( e^{\ze \MU }\) =\sum_{j}h_j{e^{\ze }\over 1-(e^{\ze }-1) ( \la_{
j}-1)}.\label{11.5u}
\end{equation}

Let $f_{ j}=1-1/\la_{ j}=(\la_{ j}-1 )/\la_{ j}$. A straightforward
calculation shows that
\begin{equation} {e^{\ze }( 1-f_{ j}) \over 1-e^{\ga }f_{ j}} ={e^{\ze }\over
1-(e^{\ze }-1) ( \la_{ j}-1)}\label{11.5v}
\end{equation} so that
\begin{equation}
\E \( e^{\ze \MU }\)=\sum_{j}h_j{e^{\ze }( 1-f_{ j})
\over 1-e^{\ze }f_{ j}}.\label{11.5w}
\end{equation}

Note that since all $\la_{ j}\geq 1/2$ we have $|f_{ j}|\leq 1$.  We can
always choose $\ze$ so that addition to (\ref{test1}) we also have
\begin{equation} |e^{\ze }|<1.\label{test2}
\end{equation} Then we can write
\begin{equation} {e^{\ze }( 1-f_{ j})
\over 1-e^{\ze }f_{ j}}=e^{\ze }( 1-f_{ j})\sum_{ k=0}^{ \ff} e^{k\ze }f^{
k}_{ j}.\label{test2.1}
\end{equation}

  Hence
\begin{equation}
\sum_{ k=1}^{ \ff}e^{k\ze } \PPP \(\MU=k\)=\E \( e^{\ze \MU }\)
=\sum_{j}h_j ( 1-f_{ j})\sum_{ k=1}^{
\ff} e^{k\ze }f^{ k-1}_{ j}.\label{11.5wa}
\end{equation}

We can choose $\ze_{ 0}<0$ so that (\ref{test1}) and (\ref{test2}) hold.
Furthermore, both sides  of (\ref{11.5wa}) are analytic functions of
$\ze$ in some neighborhood of $\ze_{ 0}+iR^{ 1}$ and agree for $\ze_{
0}+iy$ when $y$ is small. This is enough to allow us to conclude that
\be
\PPP\(\MU=k \)= \sum_{j}h_j ( 1-f_{ j})f_j^{k-1},\hspace{
.2in}k=1,2,\ldots. \label{originx}
\ee Hence
\be
\PPP\(\MU > u \)=\sum_{j}h_j f_j^{u}\hspace{ .2in}u=0,1,\ldots
\label{origin}
\ee  This completes the proof of (\ref{new:rep}).\bsq

With the aid of (\ref{new:rep}) we next provide a localization result for the
occupation measure of $\{ X_n \}$.

\begin{lemma}[The Localization Lemma]\label{lb2} Let $\{ X_n \}$  be a
symmetric transient random walk in $Z^d$ with finite second moments,
and let
$A$ be a finite set in $Z^{ d}$. Set
$\ths=\log (\La_{A}/(\La_{A}-1))$. Then for some $c_1<\infty$, $n
\geq  u^{6}$, and all
$u>0$ sufficiently large
\be c_1^{-1} e^{-\ths u} \leq
\PPP\(\mu^X_n\(A\)\geq u \)
\leq \PPP\(\mu^X_{\infty}\(A\)\geq u \)
\leq c_1 e^{-\ths u}.
\label{11.0}
\ee
\end{lemma}

\noindent {\bf Proof of Lemma \ref{lb2}:} Let $\MU_n:=\mu^X_n(A)$.
Assume first that $A$ contains the origin.  The dominant terms in
(\ref{origin}) correspond to the $f_j$'s with largest absolute value. But
since
$\PPP\(\MU > u
\)\geq 0$ and monotone decreasing it is clear that these dominant terms
must in addition be those which correspond to positive $f_j$'s. Thus the
$f_j$'s with largest absolute value are positive, i.e. correspond to
$\la_j$'s which are greater than 1. Recall that
$\phi_1$ is a strictly positive function on $A$, hence in (\ref{origin})
  we have $h_1>0$. Since $( x-1)/x=1-1/x$ is strictly monotone increasing
on
$(1,\ff)$ we conclude that the dominant term in (\ref{origin})
  is precisely the single term corresponding to the largest eigenvalue $\la_{
1}=\La_{A}$. Hence
\begin{equation}
\PPP\(\MU > u \)\sim h_1 f_1^{u}=h_1 \({\La_{A}-1 \over
\La_{A}}\)^{u}=h_1 e^{ -u\log (\La_{A} /(\La_{A}-1)  )}\label{origina}
\end{equation} implying that
\be
\lim_{u \to \infty} \PPP(\MU_\infty >u) e^{\ths u} = h_1 \in (0,\infty)
\label{eq:amir4s}
\ee out of which the upper bound of (\ref{11.0}) immediately follows.

Turning to prove the corresponding lower bound,  let $\tau_z := \inf
\{ s : |X_s|> z\}$,  and note that
\begin{equation}
\PPP(\tau_z>n)\leq c_{ 1}\exp(-c_{ 2} n z^{-2}).\label{11.simp}
\end{equation} Here is a simple proof:
\begin{eqnarray}\qquad \PPP\(\tau_z>n\) &=& \PPP\(|X_{ k}|\leq z\,;\,1\leq
k\leq n\)\nn\\ &\leq & \PPP\(|X_{ lz^{ 2}}|\leq z\,;\,1\leq l\leq nz^{-
2}\)\nn\\ &\leq & \PPP\(|X_{ lz^{ 2}}-X_{ (l-1)z^{ 2}}|\leq 2z\,;\,1\leq
l\leq nz^{-  2}\)\nn\\ &\leq &
\prod_{ l=1}^{ [nz^{- 2}]}\PPP\(|X_{ lz^{ 2}}-X_{ (l-1)z^{ 2}}|\leq
2z\)\nn\\
&=& \( \PPP\(|X_{ z^{ 2}}|\leq 2z\)    \)^{  [nz^{- 2}]}
\leq e^{ -c_{ 2} [nz^{ -2}]}.
\label{p1.2m}
\end{eqnarray} Hence
\be
\PPP(\MU_n >u) \geq \PPP(\MU_{\tau_z} >u) - \PPP(\tau_z>n) \geq
\PPP(\MU_{\tau_z} >u) - c^{-1} \exp(-c n z^{-2}) \,.
\label{eq:amir1s}
\ee

As usual we use the notation $\PPP^a$ to denote probabilities of the
random walk
$a+X_{ n},\,n=0,1,\ldots$. We now observe that
\begin{equation}
\sup_{ a\in A}\PPP^a(\mu^X_\ff (A) > u) \leq c\PPP(\mu^X_\ff (A)
>u)\label{11.assume}
\end{equation} for some $c<\ff$ and all $u$. To see this, note that for
each
$a\in A$ we can find some $n_{ a}$ with $h_{ a}=\PPP(X_{ n_{ a}}=a )>0$.
Then using the
  Markov property,
\begin{equation}
\PPP(\mu^X_\ff (A)>u)\geq \PPP(\lc\mu^X_\ff (A)>u\rc\circ\th_{  n_{
a}},\,X_{ n_{ a}}=a) =h_{ a}\PPP^{ a}(\mu^X_\ff
(A)>u).\label{11.assume3}
\end{equation}  Then (\ref{11.assume}) follows with $c=\sup_{ a\in
A}h^{ -1}_{ a}<\ff$.

Let $\MU$ and $\MU'$ denote two independent copies of $\MU_\infty$
and
$T_{A} := \inf \{s > 0: X_s \in A\}$.  Noting (\ref{11.assume}),  and  using
the strong Markov property,  it is not hard to verify that
\be
\PPP(\MU_\infty >u) \leq \PPP(\MU_{\tau_z} > u) + c\PPP(\MU+\MU'>u)
\sup_{|v|>z} \PPP^v(T_{ A}<\ff)
\label{eq:amir2s}
\ee (c.f. \cite[(3.6) and (3.7)]{DPRZ} where this is obtained for the
Brownian motion). It follows from Theorem 10.1 of \cite{Lawler2} that
\begin{equation} G( x)\leq {c\over |x|},\hspace{ .2in}|x|\neq
0.\label{Gbound}
\end{equation} Using this together with the fact that $G(X_{n\wedge
T_{A}})$ is a  martingale shows that
\begin{equation} G( v)=\E^{ v}(G(X_{ T_{ A}}),\, T_{ A}<\ff)\geq \inf_{
a\in A}G( a)
\PPP^v(T_{ A}<\ff).\label{Gmart}
\end{equation} Therefore
\begin{equation}
\sup_{|v|>z} \PPP^v(T_{ A}<\ff) \leq c z^{-1}.
\label{Gb1}
\end{equation}
  By (\ref{eq:amir4s}) it follows that  for some constant $C$ independent
of $u$, which may change from line to line,
\bea
\label{eq-csita1}
\PPP(\MU+\MU'>u)&= &
\PPP(\MU > u) +
\sum_{ y=0}^u \PPP(\MU'>u-y)
\PPP(\MU=y)\nonumber\\ &\leq & C
\Big[ \exp(-u \ths) +
\sum_{ y=0}^u \exp(-(u-y)\ths) \PPP(\MU=y)
\Big]
\nonumber \\ &\leq & C \exp(-u\ths)+C\sum_{
y=0}^u\exp(-u\ths)\nonumber \\ & =& C(1+u)\exp(-u\ths)\,.
\eea

  Hence,  taking $z= u^{2}$ one gets from (\ref{eq:amir1s}) and
(\ref{eq:amir2s}) that for some $c' >0$, all
$n \geq u^6$ and $u$ sufficiently large
\be
\PPP(\MU_{n} > u) \geq c' e^{-\ths u}
\label{eq:amir5s}
\ee as needed to complete the proof of the lemma when $A$ contains the
origin.  In general  we have
\begin{equation}
\PPP(\mu^X_\ff (A)>u)= \PPP(\lc\mu^X_\ff (A)>u\rc\circ\th_{ T_{A}},\,
T_{A}<\ff)=\sum_{a\in A}\PPP^{ a}(\mu^X_\ff
(A)>u)\PPP(T_{A}=T_{a}<\ff)\label{spread}
\end{equation} and since it is easy to see from its proof that
(\ref{eq:amir4s}) holds with
$\PPP$ replaced by $\PPP^{ a}$ for any $a\in A$, for some $ c_1= c_1(a
)$ it follows that (\ref{eq:amir4s}) also holds. This completes the proof  of
the lemma.
\bsq

{\bf Remark.}  If $A$ is replaced by $z+A$ for some fixed $z\in  Z^{ d}$,
note from (\ref{matrix}) that as matrices, $G_{z+ A}=G_{ A}$. Hence
$\La_{z+A}=\La_{A}$.

\section{Proof of Theorem 1.1}\label{sec-1.1}

Given Lemma \ref{lb2}, Theorem 1.1 follows by the methods of
\cite[Section 7]{ET}. We spell out the details.

We first prove the lower bound for (\ref{eq-slevy}). To this end fix
$a<\ths^{ -1}=1/\log (\La_{A}/(\La_{A}-1))$.

Let $k(n)=(\log n)^8$ and $N_n=[n/k(n)]$, and $t_{i,n}=ik(n) $ for
$i=0,\ldots,N_n-1$. Writing $X_s^t=X_{s+t}-X_t$ it follows that
\[
  \sup_{m \in [0,n]}
\mu^X_n(X_{ m}+A) \geq   \max_{0\leq i\leq N_n-1}\, Z_i^{(n)} \;,
\] where
$Z_i^{(n)}=\mu_{k(n)}^{X^{ t_{i,n}} }(A)$ are i.i.d. and by
Lemma \ref{lb2}, for some $c>0$ and all $n$ large enough,
\[
\PPP( \max_{0\leq i\leq N_n-1}\, Z_i^{(n)} \leq a\log n) \leq
(1-cn^{ -a \ths})^{N_n}\leq e^{-c n^{ -a \ths }N_n}\;.
\] Since $a \ths<1$ this is summable, so that applying  Borel-Cantelli, then
taking
  $a \uparrow  \ths^{ -1}$, we see that a.s.
\be
  \liminf_{ n\rar\ff}\sup_{m \in [0,n]}{ \mu^X_n(X_{ m}+A)\over  \log n}
  \geq \ths^{ -1}.\label{taylowerbound}
\ee This gives the lower bound for (\ref{eq-slevy}).

For the upper bound, fix $a> \ths^{ -1}$. Note that for any $m \in [0,n]$
\begin{eqnarray}
\mu^X_n(X_{ m}+A) &=& \sum_{ j=0}^{ n}\one_{X_{ m}+A }(X_{ j})
= \sum_{ j=0}^{ n}\one_{A }(X_{ j}-X_{ m})
\label{3.1new}\\ &=&  \sum_{ j=0}^{ m-1}\one_{A }(X_{ j}-X_{ m})+
 \sum_{ j=m}^{ n}\one_{A }(X_{ j}-X_{ m})\nonumber\\
 &\stackrel{law}{=}&  \sum_{ j=1}^{ m}\one_{A }(X'_{ j})+
 \sum_{ j=0}^{ n-m}\one_{A }(X''_{ j})\nonumber
\end{eqnarray}
where $\{ X'_{ j}\,,\,j=0,1,\ldots\}$, $\{ X''_{ j}\,,\,j=0,1,\ldots\}$ are two
independent copies of  $\{ X_{ j}\,,\,j=0,1,\ldots\}$ and we have used the
symmetry of $X_{ 1}$. Using this and (\ref{eq-csita1}),
\begin{eqnarray}
\PPP ( \sup_{m \in [0,n]}
\mu^X_n(X_{ m}+A) \geq a\log n)
  &\leq &   \sum_{ m=0}^{ n}\PPP (
\mu^X_n(X_{ m}+A) \geq a\log n)\label{3.1}\\ &\leq &
2n\PPP(\MU+\MU'\geq a\log n)\leq c(\log n)\, n^{-(a \ths  -1)
}.\nonumber
\end{eqnarray}

Thus letting $n_{ k}=n^{ k}$ for $k$ sufficiently large that $k(a
\ths  -1)>2$, we see from applying Borel-Cantelli, then taking
  $a \downarrow  \ths^{ -1}$, that a.s.
\[
  \limsup_{ n\rar\ff}\sup_{m \in [0,n^{ k}]}{ \mu^X_{ n^{ k}}(X_{
m}+A)\over
\log n^{ k}}
  \leq \ths^{ -1}.
\] The upper bound for (\ref{eq-slevy}) then follows by interpolation.

The lower bound for (\ref{seq-amir3}) follows immediately from
(\ref{taylowerbound}). As for the upper bound in (\ref{seq-amir3}),  we
note that
$\mu^X_n(x+A)=0$ unless $X_{ m}\in x+A$ for some $m \in [0,n]$. Thus
the only relevant $x$'s in (\ref{seq-amir3}) are of the form $X_{ m}-a$ for
some $m \in [0,n]$ and $a\in A$.  Thus
\be
  \sup_{x\in Z^{ d}}
\mu^X_n(x+A)= \sup_{m \in [0,n],\,a\in A}
\mu^X_n(X_{ m}-a+A).
\ee Recalling Remark and the fact that $A$ is a finite set, the  upper
bound for (\ref{seq-amir3}) now follows as in the proof of the upper
bound for (\ref{eq-slevy}).

\section{Examples}\label{sec-ex}

Proof of (\ref{ex.1}): When $A=\{ 0,y\}$ we have
\[G_{ A}=\(\begin{array}{cc} G( 0) & G( y)\\ G( y) & G(0)
\end{array}\).\] The eigenvalues are $G( 0)+G(y), G( 0)-G( y)$ so that
$\La_{A}=G( 0)+G( y)=G( 0)( 1+t_{ y})=( 1+t_{ y})/\ga_{ d}$, where
$t_{ y}=\PPP( T_{ y}<\ff)$, $\ga_{ d}$
  is the probability of no-return to the origin, and we have used the fact
that
$G( y)=t_{ y}G(0)$. Then $1-1/\La_{A}=1-\ga_{ d}/( 1+t_{ y})$.

We note that in the notation of Lemma \ref{repres}, $h_1=1,   h_2=0$ so
that by (\ref{new:rep})
\be
\PPP\(\mu^X_\infty (\{ 0,y\}) > u \)= ( 1-\ga_{ d}/( 1+t_{ y}))^{
u}\hspace{ .2in}u=1,2,\ldots \label{new:repex1}
\ee

Proof of (\ref{ex.2}): We now consider the simple random walk, and for
ease of notation consider first $d=3$. Let  $A=\{ e_{ 1},e_{ 2},e_{ 3},-e_{
1},-e_{ 2},-e_{ 3}\}=S( 0,1)$, the (Euclidean) sphere in
$Z^{ 3}$ of radius $1$ centered at the origin. We have
\[G_{ S( 0,1)}=\(\begin{array}{cccccc} G( 0) & G(e_{ 1}-e_{ 2})&G(  e_{
1}-e_{ 3}) & G(2 e_{ 1} )&G(e_{ 1}+e_{ 2})&G( e_{ 1}+e_{ 3})\\  G(
e_{2}-e_{ 1}) & G(0)&G( e_{ 2}-e_{ 3}) & G( e_{ 2}+ e_{ 1} )&G(2e_{
2})&G( e_{ 2}+e_{ 3})\\  G(e_{ 3}-e_{ 1}) & G(e_{ 3}-e_{ 2})&G( 0) & G(e_{
3}+ e_{ 1} )&G(e_{ 3}+e_{ 2})&G(  2e_{ 3})\\ G(2 e_{ 1} )&G(e_{ 1}+e_{
2})&G( e_{ 1}+e_{ 3})&G( 0) & G(e_{ 1}-e_{ 2})&G( e_{ 1}-e_{ 3})\\
  G( e_{ 2}+ e_{ 1} )&G(2e_{ 2})&G( e_{ 2}+e_{ 3}) & G( e_{2}-e_{ 1})  &
G(0)&G( e_{ 2}-e_{ 3})\\ G(e_{ 3}+ e_{ 1} )&G(e_{ 3}+e_{ 2})&G( 2e_{ 3})
&G(e_{  3}-e_{ 1}) & G(e_{ 3}-e_{ 2})&G( 0)
\end{array}\).\]

Using $G( x)=t_{ x}G(0)$, where $t_{ x}=\PPP(T_{ x}<\ff)$, and symmetry
which allows us to set $a=:t_{e_{ i}\pm e_{ j}}$  for $i\neq j$ and
$b=:t_{2e_{ i}}$ we can write

\[G_{ S( 0,1)}=G( 0)\(\begin{array}{cccccc} 1 & a&a & b&a&a\\  a & 1&a &
a&b&a\\ a & a&1 & a&a&b\\ b&a&a&1 & a&a\\ a&b&a & a & 1&a\\
a&a&b &a & a&1
\end{array}\).\]

It follows from the Perron-Frobenius Theorem that the largest eigenvalue is
$\La_{S( 0,1)}=G( 0)( 1+4a+b)$ with eigenvector
$(1,1,1,1,1,1)$. Also, it is easy to see by symmetry that $\ga_{
3}=\PPP(T_{ e_{ i}}=\ff)$.  Then again by symmetry $\PPP(T_{
e_{1}}=\ff)= { 4\over 6}\PPP(T_{e_{1}-e_{2}}=\ff)+{ 1\over
6}\PPP(T_{2e_{1}}=\ff)$, i.e.
$6\ga_{ 3}=4\PPP(T_{e_{1}-e_{2}}=\ff)+\PPP(T_{2e_{1}}=\ff)$.  Hence
$\La_{S( 0,1)}=G( 0)( 1+4a+b)=G( 0)6( 1-\ga_{ 3})=6( 1-\ga_{ 3})/\ga_{
3}$.  For the case of general $d\geq 3$, $G_{ S( 0,1)}$ is now a
$2d\times2d$ matrix, which is $G( 0)$ times a matrix in which each row
has a single entry  entry of $1$, a single entry of $b=:t_{2e_{ i}}$ and
$2d-2$ entries of $a=:t_{e_{  i}\pm e_{ j}}$, where as before $t_{
x}=\PPP(T_{ x}<\ff)$. It is easy to see by symmetry that
$\ga_{ d}=\PPP(T_{ e_{ i}}=\ff)$. Also, as before, it is easy to see  by
symmetry
$\PPP(T_{ e_{1}}=\ff)= { ( 2d-2)\over 2d}\PPP(T_{e_{1}-e_{2}}=\ff)+{
1\over 2d}\PPP(T_{2e_{1}}=\ff)$, i.e.
$2d\ga_{ d}=( 2d-2)\PPP(T_{e_{1}-e_{2}}=\ff)+\PPP(T_{2e_{1}}=\ff)$.
  Hence $\La_{S( 0,1)}=G( 0)( 1+( 2d-2)a+b)= G( 0)2d( 1-\ga_{ d})=2d(
1-\ga_{d})/\ga_{ d}$  for all $d\geq 3$.

We note that in the notation of Lemma \ref{repres}, $h_1=1,   h_j=0,
\forall j\neq 1$ so that by (\ref{new:rep})
\be
\PPP\(\mu^X_\infty (S( 0,1)) > u \)= ( 1-\ga_{ d}/2d( 1-\ga_{d}))^{
u}\hspace{ .2in}u=1,2,\ldots \label{new:repex2}
\ee  Actually, (\ref{new:rep}) assumes that $0\in A$ which doesn't  hold
here, but using (\ref{spread}) and symmetry we have that
$\PPP\(\mu^X_\infty (S( 0,1)) > u
\)=\PPP^{ e_{ 1}}\(\mu^X_\infty (S( 0,1)) > u \)$ and  (\ref{new:repex2})
follows.

Proof of (\ref{ex.3}): We again consider the simple random walk. Let  now
$A=\{0\}\cup S( 0,1)=B( 0,1)$, the (Euclidean) ball in
$Z^{ 3}$ of radius $1$ centered at the origin. With $s=\PPP( T_{ e_{
i}}<\ff)$ and
$\bar s=( s,\cdots,s)\in R^{ 2d}$ we have
\[G_{ B( 0,1)}=G( 0)\(\begin{array}{cc} 1 & \bar s\\
\bar s ^{ t}& M
\end{array}\)\] with $M$ the $2d\times 2d$ matrix in the previous
example. $M$ is a self-adjoint matrix, and as mentioned the principal
eigenvector  is $\bar 1=( 1,\cdots,1)\in R^{ 2d}$. If
$u_{ i}, i=1,\ldots, 2d-1$ denote the other orthonormal eigenvectors  of
$G( 0)M$ with  eigenvalue $\la_{ i}<\La_{S( 0,1)}$, then since they are
orthogonal to $\bar 1$ it is clear that $(0,u_{ i}), i=1,\ldots, 2d-1$ will
give us
$2d-1$ orthonormal eigenvectors of $G_{ B( 0,1)}$ with eigenvalues
$\la_{ i}<\La_{S( 0,1)}$. The remaining (two) orthonormal eigenvectors
must be of the form $( v,w\bar 1)$ and the corresponding eigenvalues will
be $G( 0)$ times those of the
$2\times 2$ matrix
\[L=\(\begin{array}{cc} 1 & 2ds\\ s  & \La
\end{array}\)\] where we abbreviate $\La=\La_{S( 0,1)}/G( 0)=2d(
1-\ga_{ d})$. The  eigenvalues of $L$ are
\begin{equation} {( 1+\La)\pm \sqrt{( 1+\La)^{ 2}-4(\La-2ds^{ 2} )} \over
2}\label{e1.0}
\end{equation} so that
\begin{equation} 1/\La_{B( 0,1)}={ 2\over G( 0)}{( 1+\La)-\sqrt{(
1+\La)^{ 2}-4(\La-2ds^{ 2} )}
\over 4(\La-2ds^{ 2} )}.
  \label{e1.1}
\end{equation} Since  $s=1-\ga_{ d}$ we have $\La-2ds^{ 2}=2d\ga_{
d}(1-\ga_{ d} )=\ga_{ d}\La$ we have
\begin{eqnarray} 1/\La_{B( 0,1)} &=& {( 1+\La)-\sqrt{( 1+\La)^{
2}-4\ga_{ d}\La}
\over 2\La}
\label{e1.3}
\end{eqnarray} so that
\begin{eqnarray}1- 1/\La_{B( 0,1)} &=& {( \La-1)+\sqrt{( 1+\La)^{
2}-4\ga_{ d}\La}
\over 2\La}
\label{e1.4a}\\&=& {( 1-1/\La)+\sqrt{( 1+1/\La)^{ 2}-4\ga_{ d}/\La}
\over 2}.\nonumber
\end{eqnarray} Setting $p=1-1/\La$ we can write this as
\begin{eqnarray}1- 1/\La_{B( 0,1)} &=& {p+\sqrt{(2-p)^{ 2}-4\ga_{
d}/\La}
\over 2}
\label{e1.4}\\&=& {p+\sqrt{p^{ 2}+4-4p-4\ga_{ d}/\La}
\over 2}
\nonumber\\&=& {p+\sqrt{p^{ 2}+2/d}\over 2}.
\nonumber
\end{eqnarray}

We note that in the notation of Lemma \ref{repres}, $ h_j=0$ for the
  $2d-1$ orthonormal eigenvectors of the form $(0,u_{ i}), i=1,\ldots,
2d-1$ above. For the principle eigenvalue we have (\ref{e1.4}) and for the
other `surviving' eigenvalue  the corresponding expression is
${p-\sqrt{p^{ 2}+2/d}\over 2}$. Hence by (\ref{new:rep})
\be
\PPP\(\mu^X_\infty (B( 0,1)) > u \)= h_{ 1}\( {p+\sqrt{p^{ 2}+2/d}\over
2}\)^{ u} +h_{ 2}\( {p-\sqrt{p^{ 2}+2/d}\over 2}\)^{ u},\hspace{
.1in}u=1,2,\ldots,
\label{new:repex3}
\ee  where $h_{ 1},h_{ 2}$ can be calculated in a straightforward manner.
We observe that since $p<\sqrt{p^{ 2}+2/d}$, the expression in
(\ref{new:repex3}) is not a mixture of geometric random variables.

Now we prove Corollary \ref{less}. For any $y\in Z^d$ we have
$t_y^2<1-\gamma_d$, since $t_y^2$ is the probability that the random walk
hits $y$ and then returns to $0$ in finite time which is obviously less
than the probability $1-\gamma_d$ that the random walk returns to zero in
finite time. This implies $(1+t_y-\gamma_d)^2<(1+t_y)^2(1-\gamma)$ which
in turn, implies
\[
-1/\log(1-\gamma_d/(1+t_y))<-2/\log(1-\gamma_d),
\]
and taking $\sup_{|y|\leq K}$ we obtain the Corollary \ref{less}.

\section{The Brownian connection}\label{sec-1.2}

Since $R_{ K}$ is a convolution operator on a bounded subset of
$\reals^d$ with locally $L^1\(\reals^d,\,dx\)$ kernel, it follows easily as
in
\cite[Corollary 12.3]{Halmos-Sunder} that $R_{ K}$ is a (symmetric)
compact operator on $L^{2}( K,\,dx)$. Moreover, the Fourier transform
relation
$\int e^{i(x\cdot p)}u^0(x)\,dx=c |p|^{-2} > 0$ implies that
$R_{ K}$ is strictly positive definite. By the standard theory for  symmetric
compact operators,
$R_{ K}$ has discrete spectrum (except near $0$) with all eigenvalues
positive, and of finite multiplicity (see \cite[Theorems VI.15,
VI.16]{Reed-Simon2}). Moreover, if we use $(f,g)_{2,K}$ to denote the
inner product in $L^{2}(  K,\,dx)$, we have that $(f,R_{ K}g)_{2,K} > 0$
for any non-negative, non-zero,
$f,g$, so by the generalized Perron-Frobenius Theorem, see
\cite[Theorem XIII.43]{Reed-Simon}, the eigenspace corresponding to the
largest eigenvalue, $\La^{ 0}_{ K}$, is one dimensional.

Let $R_{ K,\ep}$ be the operator on $L^{2}( K,\,dx)$ with kernel
\begin{equation} u^{ 0}_{ K,\ep}( x,y)=\sum_{ z\neq z'\in \LL_{ \ep}(
K)} u^{ 0}( z-z')1_{  e_{ \ep}( z)}( x)1_{  e_{ \ep}( z')}(y). \label{13.5}
\end{equation} Since the sum is over disjoint sets, it can be checked
easily that for any  $1<p<d/( d-2)$,
$u^{ 0}_{ K,\ep}( x,y)$ is bounded in $L^{ p}( K\times K,\,dx\,dy)$
uniformly
  in $\ep>0$:
\begin{eqnarray}
\int_{ K\times K} |u^{ 0}_{ K,\ep}( x,y) |^{ p}\,dx\,dy
  &=& \sum_{ z\neq z'\in \LL_{ \ep}( K)}\int_{  e_{ \ep}( z)\times e_{
\ep}( z')}
  |u^{ 0}( z-z')|^{ p}\,dx\,dy\label{13.20}\\ &=& c\ep^{ 2d}\sum_{ z\neq
z'\in
\LL_{ \ep}( K)}{ 1\over |z-z'|^{ p( d-2)}}
\nonumber\\ &=& c\ep^{ 2d-p( d-2)}\sum_{ i\neq j\in Z^{  d},\,|i|,|j|\leq
k/\ep}{ 1\over |i-j|^{ p( d-2)}}\leq C.
\nonumber
\end{eqnarray}
  Also $u^{ 0}_{ K,\ep}( x,y)\rar u^{ 0}( x-y)$ as $\ep\rar 0$ for all
$x\neq y$. Hence, using (\ref{1.edgeb})
\begin{equation}
\lim_{ \ep\rar 0}(f,R_{ K,\ep}f)_{2,K}= (f,R_{ K}f)_{2,K},\hspace{
.2in}\forall f\in C( K).\label{13.6}
\end{equation}

By \cite{Uch},
\begin{equation} G( x)=( 1+\de( x))u^{ 0}( x),\hspace{ .2in}\forall  x\neq
0\label{13.7}
\end{equation} with $\de( x)$ bounded and $\lim_{ |x|\rar\ff}\de( x)=0$
so that
\begin{equation}
\ep^{2 -d}G(\ep^{ -1} x)=( 1+\de(\ep^{ -1} x))u^{ 0}( x),\hspace{
.2in}\forall  x\in \LL_{
\ep}( K), x\neq 0.\label{13.8}
\end{equation} Let $G_{ K,\ep}$ be the operator on $L^{2}( K,\,dx)$ with
kernel
\begin{equation} v^{ 0}_{ K,\ep}( x,y)=\sum_{ z\neq z'\in \LL_{ \ep}( K)}
\ep^{2 -d} G(\ep^{ -1}( z-z'))1_{  e_{ \ep}( z)}( x)1_{  e_{ \ep}( z')}( y).
\label{13.9}
\end{equation} Using (\ref{13.8}), the same argument leading to
(\ref{13.6})
  shows that
\begin{equation}
\lim_{ \ep\rar 0}(f,G_{ K,\ep}f)_{2,K}= (f,R_{ K}f)_{2,K},\hspace{
.2in}\forall f\in C( K).\label{13.10}
\end{equation} Furthermore, since $G( 0)<\ff$, if we let $\wt G_{  K,\ep}$
be the operator on
$L^{2}( K,\,dx)$ with kernel
\begin{equation} w^{ 0}_{ K,\ep}( x,y)=\sum_{ z, z'\in \LL_{ \ep}( K)}
\ep^{2 -d} G(\ep^{ -1}( z-z'))1_{  e_{ \ep}( z)}( x)1_{  e_{ \ep}( z')}( y)
\label{13.11}
\end{equation}
  it follows from (\ref{13.10}) that
\begin{equation}
\lim_{ \ep\rar 0}(f,\wt G_{ K,\ep}f)_{2,K}= (f,R_{  K}f)_{2,K},\hspace{
.2in}\forall f\in C( K).\label{13.12}
\end{equation}

It now follows from \cite[Theorem VIII.3.6, ]{Kato} that, if $\La (
\wt G_{ K,\ep})$ denotes the largest eigenvalue of the operator $ \wt G_{
K,\ep}$
\begin{equation}
\lim_{ \ep\rar 0}\La ( \wt G_{ K,\ep})=\La (R_{ K})=\La^{
0}_K.\label{13.13}
\end{equation}

If $f$ is any eigenvector for $ \wt G_{ K,\ep}$ with eigenvalue
$\ze>0$, it is clear that $f$ is in the image of $\wt G_{ K,\ep}$ so that we
can write
\begin{equation} f( x)=\sum_{ z\in \LL_{ \ep}( K)} 1_{  e_{ \ep}( z)}( x)f(
z)\label{13.14a}
\end{equation} and the eigenvalue equation $ \wt G_{ K,\ep}f=\ze f$
becomes
\begin{equation} 1_{  e_{ \ep}( z)}( x)\sum_{ z'\in \LL_{ \ep}(  K)}\int_{
K} w^{ 0}_{ K,\ep}( x,y)1_{  e_{ \ep}( z')}( y)f( z')\,dy=\ze f( z),\hspace{
.2in}\forall z\in
\LL_{ \ep}( K).\label{13.14}
\end{equation}
  Noting that the $dy$ integration picks up a factor $\ep^{ d}$, this
implies that
\begin{equation}
\sum_{ z'\in \LL_{ \ep}( K)}
\ep^{2 } G(\ep^{ -1}( z-z'))f( z')=\ze f( z),\hspace{ .2in}\forall  z\in \LL_{
\ep}( K).\label{13.15}
\end{equation} Hence $\La_{\ep^{ -1}\LL_{\ep}( K)}= \ep^{ -2}\La ( \wt
G_{ K,\ep})$.  Together with (\ref{13.13}) this completes the proof of
Theorem
\ref{theo-conv}.

\def\noopsort#1{} \def\printfirst#1#2{#1} \def\singleletter#1{#1}
   \def\switchargs#1#2{#2#1} \def\bibsameauth{\leavevmode\vrule
height .1ex
   depth 0pt width 2.3em\relax\,}
\makeatletter \renewcommand{\@biblabel}[1]{\hfill#1.}\makeatother
\newcommand{\bysame}{\leavevmode\hbox to3em{\hrulefill}\,}

{\footnotesize

\baselineskip=12pt

\noindent
\begin{tabular}{lll} & Endre Cs\'aki
    & \hskip30pt Ant\'onia F\"oldes \\ & Alfr\'ed R\'enyi Institute of
Mathematics
    & \hskip30pt Department of Mathematics\\ & Hungarian Academy of
Sciences
    & \hskip30pt College of Staten Island, CUNY \\ & P.O. Box 127
    & \hskip30pt 2800 Victory Blvd., Staten Island\\ & H-1364 Budapest
    & \hskip30pt New York 10314 \\ & Hungary
    & \hskip30pt USA \\ & {\tt csaki@renyi.hu}
    & \hskip30pt {\tt afoldes@gc.cuny.edu} \\

& & \\
& & \\

& P\'al R\'ev\'esz
    & \hskip30pt Jay Rosen \\ & Institut f\"ur Statistik und
         Wahrscheinlichkeitstheorie
    & \hskip30pt  Department of Mathematics\\ & Technische
Universit\"at Wien
    & \hskip30pt College of Staten Island, CUNY \\ & Wiedner Hauptstrasse
8-10/107
    & \hskip30pt 2800 Victory Blvd., Staten Island \\ & A-1040 Vienna
    & \hskip30pt  New York 10314 \\ & Austria
    & \hskip30pt USA \\ & {\tt revesz@ci.tuwien.ac.at}
    & \hskip30pt {\tt jrosen3@earthlink.net} \\

& & \\
& & \\

& Zhan Shi & \\ & Laboratoire de Probabilit\'es UMR 7599 & \\ &
Universit\'e Paris VI & \\   & 4 place Jussieu & \\    & F-75252 Paris Cedex
05 & \\ & France & \\    & {\tt zhan@proba.jussieu.fr} &
\end{tabular}

}

\begin{thebibliography}{1}


\bibitem{Ciesielski-Taylor} Z.~Ciesielski and S.~J.~Taylor, {\em  First
passage times and sojourn times for Brownian motion in space and the
exact
   {H}ausdorff measure of the sample path}, Trans.\ Amer.\ Math.\ Soc.
{\bf 103}
   (1962), 434--452.


\bibitem{DPRZ} A.~Dembo, J.~Rosen, Y.~Peres and O.~Zeitouni, {\em
Thick points for spatial
   {B}rownian motion: multifractal analysis of occupation measure},  Ann.
Probab., {\bf 28} (2000), 1-35.

\bibitem{ET} P. Erd\H os and S. J. Taylor, {\em Some problems
concerning the structure of random walk paths},  Acta Sci. Hung. {\bf 11}
(1960), 137--162.
\bibitem{Halmos-Sunder} P. Halmos and V. Sunder, {\em Bounded Integral
Operators on $L^2$ Spaces}, Springer-Verlag, New York,  1978.

\bibitem{Lawler2} G. Lawler, Notes on random walks, in preparation.

www.math.cornell.edu/\verb1~1lawler/m778s04.html

\bibitem{Kato} T. Kato, {\em Perturbation Theory of Linear Operators},
second edition, Springer-Verlag, New York, 1976.

\bibitem{Reed-Simon2} M. Reed and B. Simon, {\em Methods of Modern
Mathematical Physics {I}: Functional Analysis}, Academic Press, New York,
1972.
\bibitem{Reed-Simon} M. Reed and B. Simon, {\em Methods of Modern
Mathematical Physics {IV}: Analysis of Operators}, Academic Press, New
York, 1978.

\bibitem{Uch} K. Uchiyama, {\em Green's functions for random walks on
$Z^N$}, Proc. London Math, Soc. {\bf 77} (1998), 215--240.

\end{thebibliography}
\end{document}